\documentclass[12pt]{amsart}

\def\C{\mathbb{C}}
\def\R{\mathbb{R}}

\def\Z{Z\!\!\!Z}

\def\iso{\cong}
\def\to{\longrightarrow}

\def\XxY{X\times Y}

\def\t{\mathfrak{t}}
\def\eomega{{\tilde{\omega}}}

\input xy
\xyoption{all}
\usepackage{latexsym}
\usepackage{amsmath}
\usepackage{amssymb}
\usepackage{amscd}
\usepackage{epsfig}
\usepackage{rotating}

\newtheorem{definition}{Definition}[section]

\newtheorem{ex}{Example}[section]

\newtheorem{theorem}{Theorem}[section]

\newtheorem{lemma}{Lemma}[section]
\newtheorem{lemma'}{Lemma}[theorem]

\begin{document}
\title[Reduction of Products]{Cohomology pairings on the symplectic reduction of products}
\author{R. F. Goldin} \address{Rebecca F. Goldin\\ George Mason University\\
MS3 F2, 4400 University Dr.\\
Fairfax, VA 22030\\
\texttt{rgoldin@math.gmu.edu}}
\author{S. Martin}\address{Shaun Martin\\ \texttt{shaun@mindwasabi.com}}
\thanks{R. Goldin is supported by NSF grant DMS-0305128 and by MSRI}

\begin{abstract} Let $M$ be the product of two compact Hamiltonian
$T$-spaces $X$ and $Y$. We present a formula for evaluating integrals
on the symplectic reduction of $M$ by the diagonal $T$
action. At every regular value of the moment map for $X\times Y$, the
integral is the convolution of two distributions associated to the
symplectic reductions of $X$ by $T$ and of $Y$ by $T$. Several
examples illustrate the computational strength of this
relationship. We also prove a linear analogue which can be used to find
cohomology pairings on toric orbifolds.
\end{abstract}\thanks{Math Subject Classification 53D20}

\maketitle

\section{Introduction}

Let $(M,\omega)$ be a symplectic $2n$-dimensional manifold with a Hamiltonian action by a compact torus $T$. This is equivalent to the existence of a moment map given by an invariant map
$$
\mu: M\longrightarrow \t^*
$$
satisfying
\begin{equation}\label{eq:momentmapcondition}
d\langle\mu ,x\rangle = \iota_{V_x}\omega
\end{equation}
for all $x\in \t$, where $V_x$ is the vector field on $M$ generated by
$x$. We say that $M$ is a {\em Hamiltonian $T$-space}, and we denote
$\mu$ by $\mu_M$ when the manifold is not obvious. The torus $T$ acts
locally freely on any level set $\mu^{-1}(t)$ provided $t$ is a
regular value of $\mu$. The symplectic quotient defined by
$$
M_t:=\mu^{-1}(t)/T
$$
is a manifold in the case that $T$ acts freely, and more generally may
have orbifold singularities.  $M_t$ inherits a symplectic form
$\omega_t$ from $\omega$, and so it is a symplectic orbifold. The
topology of symplectic reductions is somewhat elusive, and much work
has been centered on understanding the relationship between $M$ and
$M_t$ as $t$ varies. One avenue to understanding the topology of $M_t$
is through cobordism and wall-crossing formula (see for example
\cite{GS:birational}, \cite{Ma:cobordism}, \cite{GGK:cobordismbook},
\cite{Ka:circle}, \cite{Gr:blowup}). From the point of view of
calculating the cohomology, one can study the ring structure of
$H^*(M_t)$  (see \cite{Ki:surjectivity}, \cite{TW:abelianreductions},
\cite{Go:effective}, \cite{HK:polygonspaces} among others), or one can
study the cohomology pairings (integrals) on $M_t$ (as done in
\cite{JK:Localization}, \cite{Wi}, \cite{Gu:ReducedPhaseSpace},
\cite{JKo:symplecticcuts}). By Poincar\'e duality, both the
ring $H^*(M_t)$ and integrals on $M_t$ give the same information, but
in practice it is hard to translate from one to the other. In this
paper, we consider the integral
\begin{equation}\label{eq:integral}
\int_{M_t} \alpha e^{\omega_t}
\end{equation}
where $\alpha$ is a cohomology class on $M_t$. In the case that
$\alpha=1$, this integral is the symplectic volume of $M_t$ and is
well-known to vary polynomially in $t$
\cite{DH:reducedsymplecticform}; a formula may be found in
\cite{GLS:multiplicitydiagrams}.

We are interested in the case in which $M$ is a product of
Hamiltonian $T$-spaces, considered under the diagonal action. Many of
the objects and results we use in this paper are well-described in the
literature; our main results are an application of these theorems in
the case that $M$ is a product. This provides a tool to
compute cohomology pairings on moduli spaces such as $k$ points on $\C
P^n$, or polygon spaces. In each case, the moduli space can be realized
as a symplectic quotient of a product of Hamiltonian spaces.
The strength of our formula is that it lends itself readily to
calculation. For example, when the individual spaces are toric
varieties (as in the case of $\C P^k$'s), the product formula
describes the integral on the symplectic reduction of the product (by
the diagonal action) in terms of characteristic functions (1 or 0) on
the moment polytopes of the individual terms. Our first example will
show how this procedure works when $M$ is a product of $S^2$'s. In the
case that $M$ is a product of vector spaces, our formula provides a
means of integrating pairings on toric orbifolds using techniques
similar to those found in \cite{GS:coefficients} for toric manifolds.

While our results describe the behavior of Hamiltonian $T$ spaces,
integrals of reduction by a general Lie group $G$ may be deduced by
using results of the second author \cite{Ma:GtoT}. Suppose that $G$ is
a connected compact Lie group, and that $T$ is a fixed maximal
torus. Suppose that $G$ acts on $M$ in a Hamiltonian fashion.  Let
$\mu_G$ be a moment map for the $G$ action on $M$, and $M/\!/G:=
\mu_G^{-1}(0)/G$. Assume also that $M/\!/G$ is a manifold, as is
$M/\!/T:= M_0$, the symplectic reduction by $T$ at 0. Then
$$
\int_{M/\!/G} a = \frac{1}{|W|}\int_{M/\!/T} \tilde{a}\cup e
$$
where $W$ is the Weyl group and $\tilde{a}$ satisfies the property that its restriction to $\mu_G^{-1}(0)/T$
equals the pull-back of $a$ to the same space under the projection $\mu_G^{-1}(0)/T\rightarrow M/\!/G$. The
class $e$ is defined as follows. Let $\Delta$ be the set of roots of $G$. To each root $\alpha$ let
$L_{\alpha}:= \mu_T^{-1}(0)\times_T \C_{(\alpha)}$, where $\C_{(\alpha)}$ is the 1-dimensional representation of
$T$ given by the root $\alpha$. Then $L_\alpha$ is a line bundle over $M/\!/T$ with Euler class $e(\alpha)$.
Define $e= \prod_{\alpha\in \Delta} e(\alpha)$.

Before we state the main theorems, we recall two relevant features of equivariant cohomology. We discuss equivariant cohomology in more detail in Section \ref{se:cartanmodel}.
The inclusion $\iota:\mu^{-1}(t)\hookrightarrow M$ induces a map on
equivariant cohomology $\iota^*:H_T^*(M)\rightarrow
H_T^*(\mu^{-1}(t))$. When $T$ acts locally freely and the cohomology
is in complex coefficients, $H_T^*(\mu^{-1}(t))\iso H^*(M_t)$. The
composition of these maps is termed the {\em Kirwan map}
$$
\kappa_t:H^*_T(M)\longrightarrow H^*(M_t).
$$
and is well known to be a surjection \cite{Ki:surjectivity}.

\begin{definition} Suppose that $a\in H_T^*(M)$. We define a function $I_M(a)$ on the set
of regular values of $\mu$ on $\mathfrak{t}^*$ by
\begin{equation}\label{eq:integral2}
I_M(a)(t):= \int_{M_t}\kappa_t(a)e^{\omega_t}.
\end{equation}
\end{definition}
\noindent Thus $I_M(a)$ is the integral over the reduced space of $\kappa_t(a)$ multiplied by (a constant and)
an appropriate power of the reduced symplectic form. It should be noted that $\omega_t$ is {\em not}
$\kappa_t(\tilde{\omega})$, where $\tilde{\omega}$ is a (particular) equivariant symplectic form: see the proof
of Lemma~\ref{le:symplecticvariation}. Note also that $I_M(a)=0$ if $a$ is homogeneous of degree $k$, with $k>
\dim M_t$.

 Suppose that $\alpha$ is a homogeneous class of top degree on $M_t$, and let $a\in H_T^*(M)$ be such
that $\kappa_t(a)=\alpha$. Then Expression (\ref{eq:integral}) may be written
$$
\int_{M_t}\alpha e^{\omega_t} =  \int_{M_t} \alpha = I_M(a)(t).
$$
As $t$ varies, there is a natural diffeomorphism among $M_t$. Thus this equality holds in a neighborhood of a
regular value of $\mu$.

Let $M=X\times Y$ be a product of two compact Hamiltonian $T$ spaces $X$ and $Y$, considered as a $T$-space
under the diagonal action.  In Appendix A we give a proof of the equivariant K\"unneth theorem:
\begin{equation}\label{eq:equivariantKunneth}
H_T^*(M) = H_T^*(X)\otimes_{H_T^*(pt)} H_T^*(Y)
\end{equation}
We use $\boxtimes$ to denote this tensor product over the module $H_T^*(pt)$.  We may now state the main
results.

\vspace{.05in}

\noindent {\bf Theorem 4.1}  (A convolution formula for compact manifolds). {\em Let $X$ and $Y$ be compact
Hamiltonian $T$-spaces with isolated fixed point sets $X^T$ and $Y^T$, respectively. Let $M=X\times Y$, and
consider $M$ a Hamiltonian $T$ space under the diagonal action. For any class of the form $a\boxtimes b\in
H^*_T(M)$, there exist distributions $J_X(a), J_Y(b)$ and $J_M(a\boxtimes b)$ on $\t^*$ such that
$$
J_M(a\boxtimes b)(t)= \frac{1}{(2\pi)^k}J_X(a)\ast J_Y(b) (t)
$$
and, for regular values of $t$
$$I_M(a\boxtimes
b)(t)=\frac{1}{(2\pi)^k}J_M(a\boxtimes b)(t).$$

\vspace{.05in}

\noindent The distributions $J_X(a), J_Y(b),$ and $J_M(a\boxtimes b)$ are defined by Equation
(\ref{eq:J-distribution}).} $J_X(a)$ involves fixed point data of the $T$ action on $X$, such as the moment map
values at the fixed points, the weights of the $T$ action on the tangent space at the fixed points, and the
values of the equivariant cohomology class $a$ restricted to the fixed points. From the point of view of its
definition, it is a very different object than the integral function $I_X(a)$.

A linear version of this theorem can be deduced from results found in \cite{GLS:multiplicitydiagrams}. Let $V$
and $W$ be symplectic vector spaces with a proper moment map $\mu$ on $V\times W$. Again, we save the definition
of $J_V(a)$ for Section \ref{se:linear}.

\vspace{.05in}

\noindent {\bf Theorem 5.2}  (A convolution formula for symplectic vector spaces). {\em Let $V$ and $W$ be
symplectic vector spaces with Hamiltonian torus actions by $T$. Let $\alpha_1,\dots,\alpha_n$ be the weights of
this action on $V$ and $\beta_1,\dots,\beta_m$ be the weights on $W$. Assume that there is a proper moment map
$\mu:V\times W\rightarrow \t^*$ for the diagonal $T$ action on $V\times W$. Let $a\boxtimes b$ be an equivariant
cohomology class on $V\times W$. Then
$$
J_{V\times W}(a\boxtimes b)(t) = \frac{1}{(2\pi)^k}J_V(a)\ast J_W(b).
$$
and,  for regular values $t$ of $\mu$,
$$
I_{V\times W} (a\boxtimes b)(t) = \frac{1}{(2\pi)^k}J_{V\times W}(a\boxtimes b)(t).
$$
}

\vspace{.05in}

 The compact version may also be deduced from the linear case using cobordism, which we discuss in Section
\ref{se:cobordism}.

Each distribution $J_M$ is, up to a constant, the Fourier transform of the pushforward to a point of an
appropriate equivariant cohomology class. The localization theorem is a formula for the pushforward as a sum of
rational functions.  The main contribution of this article is to show that, in the case of product manifolds, we
gain enormous computational strength using these formulas and their Fourier transforms. The linear analogue
provides a tool for finding cohomology pairings on toric orbifolds.

\section{Background}
\subsection{Equivariant cohomology in the Cartan Model}\label{se:cartanmodel}

Let $M$ be a compact manifold of dimension $2n$ with a smooth action
by a $k$-dimensional compact torus $T$. Denote the Lie algebra of $T$
by $\mathfrak{t}$. The $T$-equivariant cohomology of $M$ may be
defined using the Cartan model: let
\begin{equation}\label{def:compactequivariantcohomology}
\Omega_T^*(M)=\Omega^*(M)^T\otimes S(\t^*)
 \end{equation}
where $S(\t^*)$ is the set of polynomials on $\t$, and
 $\Omega^*(M)^T$ is the set of smooth invariant differential forms on $M$.
  The differential on the complex may be defined by
$$
 d_T= d\otimes 1 - i\sum_{j=1}^k \iota_{V_j}\otimes u_j
 $$
 where $\{e_j\}, j=1,\dots, k$ is an orthonormal basis for $\t$ and
$u_j$ are coordinate functions such that $u_j(x)=x_j=(e_j, x)$ for any $x\in \t$, and $i=\sqrt{-1}$. The vector
fields $V_j$ on $M$ are generated by the action of $e_j$, and $\iota_{V_j}$ is an anti-derivation which
evaluates forms on $V_j$. In particular, let $\alpha\in \Omega^*(M)^T$ and $f\in S(\t^*)$. Then
$$d_T(\alpha\otimes f) = d\alpha\otimes f-i\sum_{j=1}^k
\iota_{V_j}\alpha\otimes u_jf.
$$
We assign elements of $S^l(\t^*)$ the degree $2l$. In particular, $\deg u_j = 2$, so the differential increases degree by 1.

The cohomology ring $H_T^*(M)$ is then the cohomology of the complex $\Omega_T^*(M)$ with the differential $d_T$.

Suppose that $M$ is equipped with a symplectic form
$\omega$, a Hamiltonian group action $T$, and moment map
$\mu:M\to\t^*$.
There exists a $d_T$-closed {\em equivariant} symplectic form
\begin{align*}
\eomega &= \omega+i\langle \mu, x\rangle\\
&= \omega\otimes 1 +i\sum_j\langle \mu,e_j\rangle\otimes u_j.
\end{align*}
 A direct computation and the moment map condition
 (\ref{eq:momentmapcondition})  show that $d_T\eomega=0$.

\subsection{The pushforward map}

Suppose $M$ is a compact manifold with a smooth $T$ action, and
consider the map $p:M\to pt$.
The  pushforward map
$$
p^M_*: H_T^*(M)\longrightarrow H_T^*(pt)
$$
in the Cartan model can be described explicitly using ordinary integration of ordinary forms.
Let $a=\sum_{\mathbf{j}}a_{\mathbf{j}}u^{\mathbf{j}}$ where $\mathbf{j}$ is a multi-index $\mathbf{j}=(j_1,\dots, j_k)$, $a_{\mathbf{j}}$ is a de Rham form on $M$, and $u^{\mathbf{j}}=u_1^{j_1}\cdots u_k^{j_k}$. Then the pushforward map may be defined by
\begin{equation}\label{eq:pushforward}
p^M_*(a) := \sum_{\mathbf{j}}\left(\int_M a_{\mathbf{j}}\right)u^{\mathbf{j}}.
\end{equation}
 In particular, the integral picks out the degree $2n$ component of
each $a_{\mathbf{j}}$.
The pushforward map can be extended to formal power series of
polynomials on $\t$ (elements of $S(\t^*)$) with values in
$\Omega^*(M)$ by applying the pushforward (\ref{eq:pushforward})
term-wise. Suppose that $a$ is a homogeneous form of degree
$l$. The right hand side of (\ref{eq:pushforward}) is 0 if $l<\dim
M$. If $l> \dim M$, it is possible that the pushforward is nonzero,
for the right hand side may pick up a term in $\Omega^{\dim
M}(M)\otimes S^j(\t)$, where $j>0$.

The localization formula expresses this pushforward as a sum of
integrals on fixed point sets. We state the theorem here for isolated fixed points.
\begin{theorem} \cite{AB:momentmap},\cite{BV:localisation}
Let $M$ be a compact
manifold with a smooth $T$ action, where $T$ is a compact
torus. Assume the fixed point set $M^T$ is isolated. Then
\begin{equation}\label{eq:fixedpointtheorem}
p_*^M(a) = \sum_{F\in M^T} \frac{\iota_F^*(a)}{e_F}
\end{equation}
where $\iota_F:F\hookrightarrow M$ and $\iota^*_F$ is the induced map on cohomology, and $e_F$ is the equivariant Euler class of the normal bundle to $F$ in $M$.
\end{theorem}
We briefly address the non-isolated fixed point case in Section~\ref{se:nonisolated}.

The expressions on the right hand side in (\ref{eq:fixedpointtheorem})
may be simplified. The inclusion map $\iota_F^*(a)$ picks out of
each $a_{\mathbf{j}}$ the degree-0 part
$a^{(0)}_{\mathbf{j}}$ and evaluates at $F$.
The equivariant Euler class $e_F$ is the product of the
weights $\gamma_1^F,\dots, \gamma_n^F$ (possibly with repeats)
occurring in this linear representation of $T$ on $T_FM$. Thus we rewrite
\begin{equation}\label{eq:pushformula}
p_*^M(a) = \sum_{F\in M^T}  \frac{\sum_{\mathbf{j}}a^{(0)}_{\mathbf{j}}(F)u^{\mathbf{j}}}{\prod_{l=1}^n\gamma_l^F}.
\end{equation}
In this light perhaps $p^M_*(a)$ is more obviously a function on
$\t$. The pushforward of the ``formal" form $a\wedge e^{\eomega}$ is
evaluated by expanding $e^{\eomega}$ into its
Taylor series.

\section{The Fourier transform of the pushforward map and Jeffrey-Kirwan abelian localization}\label{se:Guillemin}

We begin by establishing conventions. Let $\mathcal{S}$ denote the space of Schwartz functions. We define the Fourier transform of $\phi\in \mathcal{S}$ by
$$
\mathcal{F}(\phi)(t)=\int_{\R^k} e^{-ix\cdot t}\phi(x)dx
$$
and the inverse Fourier transform by
$$
\mathcal{F}^{-1}(\psi)(x) = \frac{1}{(2\pi)^k} \int_{\R^k} e^{ix\cdot t}\psi(t)dt.
$$
Then the Fourier transform of an arbitrary tempered distribution $h$ (continuous linear function on $\mathcal{S}$) is defined by
$$
\mathcal{F}(h)(\psi) = h(\mathcal{F}(\psi)), \hspace{.2in} \psi\in \mathcal{S}.
$$
We follow \cite{Gu:ReducedPhaseSpace} in constructing a distribution
which serves as a Fourier transform of the pushforward map, after a
choice of  {\em polarization} of the $T$ action on
$M$. Consider the distribution $h_\gamma$, with $\gamma\in \t^*$,
given by
\begin{equation}\label{eq:h-dist}
\langle h_\gamma, \varphi\rangle = \int_0^\infty \varphi(s\gamma)ds
\end{equation}
 for any (smooth, compactly supported) test function $\varphi$ on
$\t^*$. For any fixed point $F\in M^T$, let $\{\gamma_j^F\}_{j=1}^n$
be the set of weights occurring by the $T$ action on the tangent space
$T_FM$. We polarize the weights as follows. Choose $\xi\in\t$ such
that $\gamma_j^F(\xi)\neq 0$ for any $j$ and any $F$. Then let
$\gamma_j^{F,\xi} = \epsilon_j^F\gamma_j^F$, where $\epsilon_j^F = 1$
if $\gamma_j^F(\xi)> 0$ and $\epsilon_j^F = -1$ if $\gamma_j^F(\xi)<
0$. It is then clear there is a half-space in $\t^*$ on which all the
distributions $h_{\gamma_j^{F,\xi}}$ are supported. For any fixed $F$,
the convolution $h_{\gamma_1^{F,\xi}}\ast\cdots\ast
h_{\gamma_n^{F,\xi}}$ is well-defined and supported on this same half-space.

Let $\epsilon^F= \prod_j \epsilon_j$. Let $a = \sum a_{\mathbf{j}}
\otimes u^{\mathbf{j}}$ be an equivariant class on $M$, where ${\bf
j}$ is a multi-index, and $a_{\mathbf{j}}$ are differential forms on $M$. For
${\bf j}=(j_1,j_2,\dots,j_k)$, let $\partial^{\bf j}$ be the differential
operator given by
\begin{equation*}
\partial^{\bf j} = \left(\frac{\partial}{\partial
u_1}\right)^{j_1}\left(\frac{\partial}{\partial
u_2}\right)^{j_2}\cdots \left(\frac{\partial}{\partial
u_k}\right)^{j_k}.
\end{equation*}
Finally, we define
\begin{equation}\label{eq:J-distribution}
J_M(a) = (2\pi)^k\sum_{F\in M^T} \epsilon^F  \left(\sum_{\bf j} a^{(0)}_{\bf
j}(F)(-i)^{|{\bf j}|}\partial^{\bf j}\right)(\delta_{\mu(F)}\ast
h_{\gamma_1^{F,\xi}}\ast\cdots\ast h_{\gamma_n^{F,\xi}}),
\end{equation}
where
$\delta_{\mu(F)}$ is the delta distribution centered at $\mu(F)$, and
$|{\bf j}|=j_1+\dots+j_k$. The following results set the context for
this paper. The first is a version of the
Jeffrey-Kirwan Abelian localization theorem \cite{JK:Localization}.

\begin{theorem}[Jeffrey-Kirwan]\label{th:JeffreyKirwan} Let $M$ be a compact Hamiltonian $T$ space with isolated fixed points. Let $I_M(a)$ be defined as in Equation (\ref{eq:integral}). If $t$ is a regular value of the moment map, then
$$
I_M(a)(t) =\frac{1}{(2\pi)^k}J_M(a)(t).
$$
\end{theorem}

The following result can be found in \cite{Gu:ReducedPhaseSpace} and is what makes the distributions $J_M(a)$ so readily computable.

\begin{theorem}[Guillemin]\label{th:Guillemin}
Let $M$ be a compact symplectic manifold. Suppose $T$ acts on $M$ in a
Hamiltonian fashion with isolated fixed points. Let $f_a(x)=
p_*^M(a\wedge e^{\eomega})(x)$ be the pushforward of the formal
equivariant form $a\wedge e^{\eomega}$ to a point. Then
$f_a$ is a function of
tempered growth whose Fourier transform is
$$
\mathcal{F}(f_a)(t)=(2\pi i)^nJ_M(a)(t).
$$
In particular, $J_M(a)$ is independent of choice of
polarization.\footnote{Note that individual terms of the ABBV fixed
  point theorem may not so easily be interpreted in this
  light. However the
  rational term $\frac{1}{\gamma}$ may be viewed as the generalized
  function $\frac{1}{\gamma+0i}$ as in \cite{GeS}.}
\end{theorem}

\section{The reduction of a product of compact Hamiltonian $T$-spaces}

We now restrict our attention to the case that
$$
M= X\times Y
$$
where $X$ and $Y$ are Hamiltonian $T$-spaces with isolated fixed
points, and $M$ is considered under the diagonal action on the
product. Our main result in the compact case is Theorem \ref{th:maincompact}, which
may be viewed as a corollary to Theorem
\ref{th:Guillemin}, or to the cobordism results of
\cite{GGK:cobordismbook}.In section
\ref{se:examples}
we present two calculations using Theorem \ref{th:maincompact}.

\subsection{The pushforward satisfies a product rule}

The class of the equivariant symplectic form on $M$ may be written
$$[\eomega_M] = [\eomega_X]\boxtimes 1 + 1\boxtimes [\eomega_Y]$$
where $\eomega_X$ and $\eomega_Y$ are equivariant symplectic forms on
$X$ and $Y$, respectively and the tensor product is that in (\ref{eq:equivariantKunneth}).
For simplicity, we assume that
$\alpha=a\boxtimes b$, and
consider the pushforward of $\alpha \cdot e^{[\eomega_M]}$ to a point.

The fixed points of $X\times Y$ are pairs $(p,q)$, where $p\in X^T$ and $q\in Y^T$. By (\ref{eq:fixedpointtheorem}) we have
\begin{align*}
p^{M}_*(\alpha e^{\eomega_M})&=p^{X\times Y}_*(a\boxtimes b e^{\eomega_X\boxtimes 1 + 1\boxtimes \eomega_Y})\\
                                &=p^{X\times
                                Y}_*(ae^{\eomega_X}\boxtimes
                                be^{\eomega_Y})\\
                                &=\sum_{(p,q)\in X^T\times Y^T}
                                \frac{\iota_{(p,q)}^*(ae^{\eomega_X}\boxtimes
                                be^{\eomega_Y})}{e_{(p,q)}}
\end{align*}
where $e_{(p,q)}$ is the equivariant Euler class of the normal bundle
to $(p,q)$ in $X\times Y$. Since
$\nu_{(p,q)}(X\times Y)=\nu_p(X)\oplus \nu_q(Y)$, we have $e_{(p,q)} =
e_p\cdot e_q$, where $e_p$ and $e_q$ are the equivariant Euler classes
of $\nu_p(X)$ and $\nu_q(Y)$, respectively. To simplify the numerator, we
expand out the terms at the level of forms and using the formula
(\ref{eq:fixedpointtheorem}) to obtain
$$
\iota_{(p,q)}^*(ae^{\eomega_X}\boxtimes be^{\eomega_Y})= \iota_p^*(ae^{\eomega_X})\iota_q^*(be^{\eomega_Y}).
$$
We have shown:
\begin{lemma} For all cohomology classes $a\in H_T^*(X)$ and $b\in H_T^*(Y)$,
\begin{equation}\label{eq:productpushforward}
p_*^{X\times Y}(a\boxtimes be^{\eomega_{X\times Y}})= p_*^X(ae^{{\eomega_{X}}})p_*^Y(ae^{{\eomega_{Y}}}).
\end{equation}
\end{lemma}

Choose a polarization for the weights at all the fixed points of
$M$. The weights appearing in the product are the union of the weights
of the action on the tangent spaces to fixed points of $X$ and of
those to fixed points of $Y$. Thus a polarization of $M$ is also a
polarization of $X$ and of $Y$. We Fourier transform Equation
(\ref{eq:productpushforward}). By construction the
distributions $h_\gamma$ associated to the polarized weights of $X$
and those of $Y$ are supported on the same convex set. Thus the right hand side
is a convolution of the Fourier transform of each of the two
pushforward maps. In particular, we have shown that
\begin{theorem}\label{th:maincompact}  Let $X$ and $Y$ be compact Hamiltonian $T$-spaces
with isolated fixed point sets $X^T$ and $Y^T$, respectively. Let
$M=X\times Y$, and consider $M$ a Hamiltonian $T$ space under the
diagonal action. Let $J_X$, $J_Y$ and $J_M$ be defined as in Equation (\ref{eq:J-distribution}). For any class of the form $a\boxtimes b\in H_T^*(M)$,
$$
J_M(a\boxtimes b)(t)= \frac{1}{(2\pi)^k}J_X(a)\ast J_Y(b) (t)
$$
and, for regular values of $t$,
$$J_M(a\boxtimes b)(t)= (2\pi)^kI_M(a\boxtimes b)(t).$$
\end{theorem}

\section{The linear case}\label{se:linear}
Suppose now that $V$ is a $2n$-dimensional symplectic vector space
with Hamiltonian $T$-action and proper moment map
$$
\mu: V\longrightarrow \t^*.
$$ We assume also that 0 is the only fixed point in $V$. The image
$\mu(V)$ lies in the convex cone given by $\{\mu(0)+\gamma_j\}$ where $\{\gamma_j\}$ are the weights of the $T$
action on $V$. Suppose that $t$ is a regular value of $\mu$ and as before let $V_t:= \mu^{-1}(t)/T$ be the
symplectic reduction of $V$ at $t$. The spaces $V_t$ are termed weighted projective varieties when $T\cong S^1$
and toric orbifolds more generally. The Kirwan map
$$
\kappa_t:H_T^*(V)\longrightarrow H^*(V_t)
$$
remains surjective since $\mu$ is proper. Since $V$ is contractible,
$$
H_T^*(V)\iso H_T^*(pt)\iso \C[u_1,\dots, u_k]$$
 where $k=\dim T$. Then any class on $V_t$ is of the form $\kappa_t(a)$ where $a=\sum a_{\bf j}u^{\bf j}$ is a polynomial in $u_1,\dots, u_k$ with complex coefficients.

The linear analogue of Theorem \ref{th:maincompact} is not independent of polarization, although the integrals
one wants to compute are. Although the existence of such a formula may seem intuitive given the presentation of
the distribution $J_M(a)$, the proof of Guillemin's theorem (Theorem \ref{th:Guillemin}) relies on the fact that
$M$ is compact and, by the Jeffrey-Kirwan theorem, $J_M(a)$ is compactly supported. In the linear case, the
relevant distributions are supported on cones, and the Jeffrey-Kirwan theorem does not hold. In addition, there
is no well-defined push-forward map for non-compactly supported equivariant forms.

The key point in the linear case is that every cohomology class on $V_t$ is obtained by an appropriate
derivative of $e^{\omega_t}$, where $\omega_t$ is the inherited symplectic volume form on $V_t$. Indeed, this is
the tactic adopted in \cite{GS:coefficients}, where the authors use these derivatives as a new technique to find
the cohomology rings of (smooth) toric varieties, following the original work by \cite{Da:toricvarieties}. We
begin by proving the following lemma.

\begin{lemma}\label{le:symplecticvariation}
The forms $\omega_t$
depends linearly on $t$ in a connected region of regular values of $\mu$. More specifically, $\frac{\partial}{\partial t_\beta} \omega_t = -i\kappa_t(u_\beta)$
\end{lemma}
\begin{proof} Let $\eomega_t = \omega\otimes 1 +i\sum_{j}(\langle
\mu,e_j\rangle-t_j)\otimes u_j$, where $t_j$ are the constant
coordinate functions on $M$. A direct calculation shows that
$d_T\eomega_t=0$, and that $\kappa_t(\eomega_t)=\omega_t$ for regular
values of $\mu$.
Note that $\eomega_t = \eomega-i\sum t_j\otimes u_j$, where
$\eomega = \omega\otimes 1 + i\sum\langle \mu,e_j\rangle\otimes
u_j$. Thus $\eomega = \eomega_t+i\sum t_j\otimes u_j$ and
$$
\kappa_t(\eomega)=\omega_t+i\sum_j t_j\kappa_t(u_j)=\omega'.
$$ Here $\omega'$ is some locally constant class over $t$ on $M_t$, as
$\kappa_t$ is a map on cohomology. Thus $\omega_t=\omega'-i\sum
t_j\kappa_t (u_j)$ and
$
\frac{\partial}{\partial t_\beta} \omega_t = -i\kappa_t(u_\beta).
$
\end{proof}
Since integration over $V_t$ and the Kirwan map $\kappa_t$ are
topological maps, we obtain
$$\frac{\partial}{\partial
t_\beta}\int_{V_t}e^{\omega_t}=
\int_{M_t}\frac{\partial}{\partial t_\beta}e^{\omega_t} =
-i\int_{M_t}e^{\omega_t}\kappa_t(u_\beta). $$
Note that this argument works as well for a compact manifold $M$.\footnote{Differentiating $(\dim M_t)$ times, one obtains the integral of a class
independent of $t$, which is constant.
It follows $\int_{M_t}e^{\omega_t}$ is locally polynomial in $t$ with degree less than or equal to $\dim M_t$. This is the content of the Duistermaat-Heckman theorem \cite{DH:reducedsymplecticform}.}

It follows that for any monomial $a= c\cdot u^{\bf j}$ with $c\in \C$
and ${\bf j}$ a multi-index,
\begin{align*}
\kappa_t(a) &= c (i)^{|\bf j|}\partial^{\bf j}\frac{(\omega_t)^{|{\bf j}|}}{|{\bf j}|!}\\
\end{align*}
where $|{\bf j}| = j_1+\dots +j_k$. Then
$$
\int_{V_t} \kappa_t(ae^{\eomega_t})=\int_{V_t}\kappa_t(a)e^{\omega_t} = \int_{V_t}c(i)^{|{\bf j}|}\partial^{\bf j}e^{\omega_t}
$$
on any cone of regular values of $\mu$. We now analyze this integral in light of what is understood about symplectic volume.


 The condition that $\mu$
be proper assures us that the weights $\gamma_1,\dots, \gamma_n$ of
the $T$ action on $V$ are polarized, i.e. there exists $\xi\in \t$
such that $\gamma_i(\xi)>0$ for $i=1,\dots, n$. Define $h_{\gamma}$ as
in Expression (\ref{eq:h-dist}). In \cite{GLS:multiplicitydiagrams}, the authors
show that
\begin{equation}\label{eq:sympvol}
(2\pi)^k\delta_{\mu(0)}\ast h_{\gamma_1}\ast\cdots\ast h_{\gamma_n}
\end{equation}
is a distribution equal to the symplectic volume of $V_t$ in the case
that the $T$-action is free. In the locally free case, at least one
vector $\gamma_i=m_i\alpha_i$, where $m_i\in \Z$ is greater than 1 and
$\alpha_i$ is a unit weight. In this case,
$h_{\gamma_i}=\frac{1}{m_i}h_{\alpha_i}$ (which can be checked by direct calculation) so
the term (\ref{eq:sympvol}) automatically picks up a term
$\frac{1}{|S|}$, where $S$ is the generic stabilizer of $T$ on $V$. We note also that the convolution is well-defined because there is a half-space on which all $h_{\gamma_j}$ are supported. By differentiating the distribution (\ref{eq:sympvol}) we obtain
$$
(-i)^{|{\bf j}|}\partial^{\bf j}\delta_{\mu(0)}\ast h_{\gamma_1}\ast\cdots\ast h_{\gamma_n}(t) = \int_{V_t} \kappa_t(u^{\bf j})e^{\omega_t}.
$$
We define a distribution on $\t^*$ by
\begin{equation}\label{eq:vectorspaceterm}
J_V(a)= (2\pi)^k \sum_{\bf j} a_{\bf j}(-i)^{\bf |j|}D^{\bf j}(\delta_{\mu(0)}\ast h_{\gamma_1}\ast\cdots\ast h_{\gamma_n})
\end{equation}
where $a=\sum_{\bf j} a_{\bf j}u^{\bf j}$ is an equivariant cohomology class on $V$. We obtain
the following linear version of the Jeffrey-Kirwan Abelian localization theorem:
\begin{theorem}
Let $V$ be a $2n$-dimensional symplectic vector space with Hamiltonian $T$ action, an isolated fixed point at 0, and a proper moment map $\mu$. Let $\gamma_1,\dots, \gamma_n$ denote the weights of the action, and let $V_t$ denote the symplectic reduced space for $t$ a regular value of $\mu$. Then for all such regular values,
$$
J_V(a)(t) = (2\pi)^k\int_{V_t}\kappa_t(a)e^{\omega_t}
$$
for any equivariant cohomology class $a\in H_T^*(V)$.
\end{theorem}
Note that (\ref{eq:vectorspaceterm}) looks like an
individual term in (\ref{eq:J-distribution}), but while the individual
terms of (\ref{eq:J-distribution}) depend on a choice of
polarization, the expression here does not.


We now show the linear analogue of the compact convolution formula (Theorem \ref{th:maincompact}). A new
subtlety arises without the compact
assumption. Let $(V, \omega_V)$ and
$(W,\omega_W)$ be two symplectic vector spaces with Hamiltonian
$T$-actions and let $\mu_V$ and $\mu_W$ be moment maps for $V$ and
$W$, respectively. The diagonal $T$-action on $V$ and $W$ has a moment
map $\mu = \mu_V+\mu_W$ and may not be proper. In fact, if
$\alpha_1,\dots, \alpha_n$ are the weights of the $T$ action on $V$
and $\beta_1,\dots, \beta_m$ are those on $W$, then the convolution
$h_{\alpha_1} \ast\cdots\ast h_{\alpha_n}\ast h_{\beta_1} \cdots\ast
h_{\beta_m}$ may not be well-defined.

This can be fixed by changing the symplectic form on $W$, providing
that the moment map $\mu_V$ is proper. Let $\xi\in \t$ be such that
$\langle \alpha_j,\xi\rangle >0$ for all $j$. We note that
$W=\oplus_j\C_{\beta_j}$ where $\C_{\beta_j}$ is a 1-dimensional
(complex) representation of $T$ with weight $\beta_j$. This splitting
is symplectic, so that we may find local coordinates $(x_j,y_j)$ on
$\C_{\beta_j}$ such that
$\omega_W=\sum_jdx_j\wedge dy_j.$ Then let
\begin{equation}\label{eq:changeomega}
\omega_W' = \sum_j sgn\langle \beta_j,\xi\rangle dx_j\wedge dy_j.
\end{equation}
It follows that $\omega_{V\times W}' = \omega_V\otimes 1 +1\otimes
\omega_W'$ is a symplectic form on $V\times W$ with proper moment map
given by
$$
\mu_{V\times W}' = \mu_V+\mu_W'.
$$
where we fix $\mu_V(0)=\mu_W'(0)=0$. We may conclude
\begin{theorem}\label{th:mainlinear}
Let $V$ and $W$ be symplectic vector spaces with Hamiltonian torus actions by $T$. Let $\alpha_1,\dots,\alpha_n$
be the weights of this action on $V$ and $\beta_1,\dots,\beta_m$ be the weights on $W$. Assume also (by changing
the symplectic forms if necessary) that there is a proper moment map $\mu:V\times W\rightarrow \t^*$ for the
diagonal $T$ action on $V\times W$. Let $a\boxtimes b$ be an equivariant cohomology class on $V\times W$ that is
homogeneous of degree $\dim (V\times W)_t$. Then for regular values $t$ of $\mu$,
$$
(2\pi)^k\int_{(V\times W)_t}\kappa_t(a\boxtimes b) = J_{V\times W}(a\boxtimes b)(t)
$$
and
$$
J_{V\times W}(a\boxtimes b)(t) = \frac{1}{(2\pi)^k}J_V(a)\ast J_W(b)
$$
where the distributions $J_V(a)$ and $J_W(b)$ are given by (\ref{eq:vectorspaceterm}).
\end{theorem}
Note that the integral $\int_{(V\times W)_t}\kappa_t(a\boxtimes b)$ equals $I_{V\times W}(a\boxtimes b)(t)$ (as
defined in (\ref{eq:integral})) since $\mu$ is proper and the degree of $a\boxtimes b$ equals the dimension of
the reduced space.

While the linear case seems extremely similar to the compact case, there are some differences that should be
emphasized here. First, for the compact case, the polarization must be for all weights at all fixed points all
at once. The corresponding distribution $J_M$ is independent of the polarization, while the individual terms are
not. In contrast, in the linear case, the distribution $J_V$ (or, rather, $J_{V\times W}$) does indeed depend on
the polarization. However, for $V_t$ to have meaning, we assumed that the moment map $\mu_V$ is proper. The
moment map can be made to be proper by changing the symplectic form as described above, or, equivalently,
changing the signs of the weights of the action so that the convolution $h_{\alpha_1} \ast\cdots\ast
h_{\alpha_n}\ast h_{\beta_1} \cdots\ast h_{\beta_m}$ is well-defined. Equivalently, one can choose a
polarization as in Section 3, and orient the weights according to the polarization so that the convolution is
well-defined. Thus for linear case (unlike the compact case), we choose a symplectic structure, thereby
assigning meaning to $(V\times W)_t$ when the moment map $\mu_{V\times W}$ is not {\em a priori} proper.






\section{Cobordism}\label{se:cobordism}

Clearly the linear case provides a one-term version of the compact
case. Off the set $\gamma_j=0$ for any $j$, the inverse Fourier
transform of $J_V(a)$ is (up to a constant) the rational function
$$
\frac{a}{\prod_{l=1}^n\gamma_l} = \frac{\sum_{\bf j}a_{\bf j}u^{\bf j}}{\prod_{l=1}^n\gamma_l}
$$
which resembles a term in the localization theorem (\ref{eq:pushformula}) and lends computational strength to Theorem~\ref{th:maincompact}. However the compact case
does not follow immediately from the linear case. Indeed, the cobordism results of \cite{GGK:cobordism} combined with the results of Section~\ref{se:linear} imply the compact convolution formula.

 Let $M$ be a compact, symplectic manifold with Hamiltonian $T$ action, isolated fixed points $M^T$, and moment map $\mu: M\rightarrow \t^*$. For each fixed point $F\in M^T$, there is a linear isotropy action of $T$ on $T_FM$, and a symplectic form $\omega_F$ given by (\ref{eq:changeomega}), possibly different from the form $\omega$ on $M$ restricted to $F$. We use the same vector $\xi$ in choosing its symplectic form for each fixed point $F$. Each tangent space $T_FM$ has a proper moment map $\mu_F$ which is unique up to a choice of constant.  Fix $\mu_F(0)=\mu(F)$.

\begin{theorem}[\cite{GGK:cobordism}]\label{th:cobordism} Let $M$ be a compact, symplectic manifold with Hamiltonian $T$ action, isolated fixed points, and moment map $\mu: M\rightarrow \t^*$. Let $t$ be a regular value of $\mu$. The reduced space $M_t$ is cobordant as a symplectic orbifold to a disjoint union of compact symplectic toric orbifolds. These orbifolds are the disjoint union of $M^F_t:= \mu_F^{-1}(t)/T$ as $F$ varies over all points in $M^T$.
\end{theorem}

It should be noted that cobordism of symplectic orbifolds carries the information of the moment map on $M$ and that on each linear space $T_FM$. In general some of the reduced space $M^F_t$ are ``obviously" empty, as the point of reduction is not in the image of $\mu_F$.

Integration is a cobordism invariant, which then implies that an integral over $M_t$ can be evaluated by an
integral over the disjoint union of $M^F_t$. These latter spaces are orbifolds obtained via symplectic reduction
of a vector space (and hence the methods of Section~\ref{se:linear} apply). In the previous section we showed
how the integral of $ae^{\eomega}$ over $M_t^F$ (copies of $V_t$) are obtained from the distribution
$J_{M^F}(a)$, where the symplectic form may have been modified to ensure that this distribution is well-defined.
The way that these forms are modified may be chosen to be according to the same polarization. Thus the linear
case of the Jeffrey-Kirwan theorem and the cobordism result imply Guillemin's theorem \ref{th:Guillemin}.

In the interest of reproving the product convolution theorem \ref{th:maincompact}, we note that the fixed point
sets of the product $M=X\times Y$ are pairs of points $(p,q)$ with $p\in X^T$ and $q\in Y^T$, and the tangent
spaces satisfy $T_{(p,q)}M=T_pX\oplus T_qY$. It then follows that the integral of $\kappa_t(ae^{\eomega})$ over
$M_t$ can be obtained by sums of the same map over the reduced spaces associated to $T_pX\oplus T_qY$ with the
changed symplectic structure described. The linear case proven in Section~\ref{se:linear} shows that integrals
over the reduction of $T_{(p,q)}M$ are obtained via a convolution, and then cobordism again (along with
linearity of the convolution) lead us to \ref{th:maincompact}.

\section{Non-isolated fixed points}\label{se:nonisolated}

All of the theorems we mention here have generalizations to non-isolated fixed points, although few of the
simplifications (such as a way of writing an Euler class as a product of weights) carry over to the case of
non-isolated fixed points. First we note that Theorem \ref{th:cobordism} is true for non-isolated fixed points,
with the toric orbifolds $M^F_t$ replaced by toric orbifold bundles over each connected component $F$ of the
fixed point set. The ABBV formula for the pushforward also holds more generally, with an integral over the fixed
point set inserted for each component of the fixed point set. Thus the cobordism result combined with a linear
theorem regarding vector bundles (rather than vector spaces) would imply a version of Theorem \ref{th:Guillemin}
in the compact case with non-isolated fixed points.  Following the same proof as in this article, a non-isolated
fixed point version of \ref{th:maincompact} would then follow. For $M=X\times Y$, each component of the fixed
point set $M^T$ is a product of two connected components $P$ and $Q$ of the fixed point sets of $X$ and $Y$,
respectively. Thus it remains true that $\nu_{(P,Q)}M = \nu_PX\oplus \nu_QY$, and hence there is a
multiplicative property of their equivariant Euler classes. The key problem (still unanswered as far as the
authors are aware) then, is to prove a non-compact version of the Jeffrey-Kirwan residue formula for vector
bundles over non-isolated fixed point sets.

\section{Examples}\label{se:examples}
We explore two examples, mainly the symplectic reduction of a product of 2-spheres (by a circle) and the
reduction of a product of $\C P^2$ (by a 2-torus). The ring structures of these products have been studied by
\cite{HK:polygonspaces} in the first case, more generally by \cite{Mo:circle} for reductions of $S^1$ actions on
the product of $\C P^n$s, and by the author \cite{Go:weightvarieties} for a product of $\C P^n$s, reduced by a
nonabelian group.
\subsection{Symplectic reduction of a product of 2-spheres}

Let $(S^2,\omega)$ be the 2-sphere with symplectic form $\omega$. Let
$S^1$ act on $S^2$ in a Hamiltonian fashion with fixed point set
$\{N,S\}$. We choose $\omega$ so that there is a moment map
$$
\mu_{S^2}: S^2\longrightarrow \R
$$
with $\mu_{S^2}(N)=1$ and $\mu_{S^2}(S)=-1$.
We assume also that $S^1$ acts on the tangent space $T_NS^2$ with
weight 1 and on $T_SS^2$ with weight -1.

Recall that for any Hamiltonian $T$-space the map
$$i^*: H_T^*(M)\longrightarrow H_T^*(M^T)$$ induced by the inclusion
$M^T\hookrightarrow M$ of the fixed point set is an injection. Thus
when $M=S^2$, the equivariant symplectic class is determined by a
class on $N$ and a class on $S$. Note that in the $S^1$ case,
$$ H_{S^1}^*(N) = H_{S^1}^*(S) = \C[u]
$$
Let $\nu$ be such the equivariant symplectic class with the
properties
$$i_N^*(\nu) = u \hspace{.1in} i_S^*(\nu)=-u.$$ Note that $\nu$ is
degree 2.

We now apply Theorem \ref{th:maincompact} to $M=S^2\times\cdots\times
S^2=(S^2)^n$ with $n$ odd, and
$a=\nu^{k_1}\boxtimes\cdots\boxtimes\nu^{k_n}$, where $\sum k_j =
n-1=\frac{1}{2}\dim (S^2)^n_0$ and the tensor product $\boxtimes$ is
taken over the module $H_T^*(pt)$ (see Appendix~\ref{se:Kunneth}). We
then use this calculation to find the symplectic volume of the
reduction $M_0=(S^2)^n_0$. A similar calculation on this moduli space can be found in \cite{Ma:cobordism}. Notice that $n$ odd implies that 0 is a regular
value of $\mu_M = \mu_{S^2}+\cdots +\mu_{S^2}.$

Since 0 is regular,
$$
\int_{(M_0}\kappa_0(\nu^{k_1}\boxtimes \cdots\boxtimes \nu^{k_n}) =\frac{1}{2\pi}
J_{M}(\nu^{k_1}\boxtimes\cdots\boxtimes \nu^{k_n})(0)
$$
and by Theorem \ref{th:maincompact} we have
$$
J_{M}(\nu^{k_1}\boxtimes\cdots\boxtimes \nu^{k_n}) = \frac{1}{(2\pi)^{n-1}}J_{S^2}(\nu^{k_1})\ast\cdots\ast J_{S^2}(\nu^{k_n}).
$$
We evaluate this convolution: Theorem \ref{th:Guillemin} implies that
\begin{align*}
p_*^{S^2}(\nu^{k_1}e^{\eomega_{S^2}}) &=2\pi i\mathcal{F}^{-1}
(J_{S^2}(\nu^{k_1})) \\
&= \frac{u^{k_1}e^{iu}}{u} + \frac{(-u)^{k_1}e^{-iu}}{-u}\\
&= u^{k_1-1}e^{iu} + (-1)^{k_1-1}u^{k_1-1}e^{-iu}\\
&= u^{k_1-1}(e^{iu}+(-1)^{k_1-1}e^{-iu})
\end{align*}
Using the product rule for the pushfoward we find
\begin{align*}
p^{M}_*(\nu^{k_1}\boxtimes \cdots\boxtimes
\nu^{k_n}e^{\eomega_{M}})
&=  \prod_{k_j} u^{k_j-1}(e^{iu}+(-1)^{k_j-1}e^{-iu})\\
&= u^{\sum k_j
- n}\prod_m (e^{iu} + e^{-iu})\prod_{n-m} (e^{iu}-e^{-iu})\\
\end{align*}
where $m$ is the number of odd $k_j$'s and $n-m$ is the number of even
$k_j$'s. Using $\sum k_j = n-1$ we simplify this expression
\begin{align*}
&= u^{-1}\left(\sum_s {m\choose{s}} e^{isu}e^{-i(m-s)u}\right)\left(\sum_r
{n-m\choose{r}} e^{iru}e^{-i(n-m-r)u}(-1)^{n-m-r}\right)\\
& = u^{-1}\sum_{s,r} {{m}\choose{s}}{{n-m}\choose{r}}(-1)^{n-m-r}e^{i(2s+2r-n)u}.
\end{align*}

For regular values of $\mu_M$ we note that
$\mathcal{F}(p_*^{M}(ae^{\eomega_M}))(t) = 2\pi i J_M(a)(t) = 2\pi i\cdot
2\pi I_M(a)(t)$, where $a = \nu^{k_1}\boxtimes \cdots\boxtimes
\nu^{k_n}$. We Fourier transform this expression by noting
that
$\mathcal{F}\left(\frac{e^{i(2s+2r-n)u}}{2\pi i u}\right) = 2\pi
\delta_{2s+2r-n}\ast h$ up to a distribution supported on the lattice
points of $\mathfrak{t}^*$, where $h(t)$ is the step function $h(t)=0$
for $t<0$ and $h(t)=1$ for $t\geq 0$.  For $t=0$ there is a contribution
whenever $s+r<\frac{n}{2}$. Note also that if $m$ is the number of odd
$k_i$, and $n$ is odd, then $\sum k_i = n-1$ is even and $m$ must be
even. Then $n-m$ must be odd and $(-1)^{n-m-r}=(-1)^{r-1}$. Therefore,
\begin{equation}\label{eq:terminsympvolume}
\int_{(S^2)^n_0}\kappa_0(\nu^{k_1}\boxtimes \cdots\boxtimes \nu^{k_n}) = \sum_{s+r<\frac{n}{2}}{m\choose{s}}{n-m\choose{r}}(-1)^{r-1}.
\end{equation}
To find the symplectic volume of this product, we note that the equivariant
symplectic form on $M = (S^2)^n$ is
$$
\omega_M = \sum_j 1\otimes\cdots\otimes i\nu\otimes\cdots\otimes 1
$$
where $\nu$ occurs in the $j$th position of the $j$th term of the sum. The symplectic volume form on the reduced
space $(S^2)^n_0$ is obtained from the Kirwan map $\kappa_0$ applied to
$\frac{(\omega_M)^{n-1}}{(n-1)!}$.\footnote{Note that this statement is only true in the case of reduction at
0.} Thus the symplectic volume is a sum of terms of the form (\ref{eq:terminsympvolume}).

\subsection{Symplectic reduction of a product of $\C P^2$'s}
Consider now the symplectic manifold $(\C P^2,\omega)$ with a $T^2$ Hamiltonian
action fixing points $\{N,S,E\}$ and moment map
$$
\mu_{\C P^2}: \C P^2\longrightarrow \R^2
$$
in which $\mu_{\C P^2}(N) = (-1,2)$, $\mu_{\C P^2}(S) = (-1,-1)$, and $\mu_{\C P^2}(E) = (2,-1)$. We use Theorem
\ref{th:maincompact} to find the symplectic volume of the reduced space $(\C P^2\times \C P^2)_0$.

We describe the equivariant symplectic form on $\C P^2$ by its
restriction to the fixed point set, where
$$
H_T^*(N) = H_T^*(S) = H_T^*(E) = \C[u_1,u_2].
$$
Choose $\nu$ with the properties that
$$\
\iota_N^*(\nu) = -u_1+2u_2, \hspace{.1in} \iota_S^*(\nu) = -u_1-u_2,
\mbox{ and }\iota_E^*(\nu) = 2u_1-u_2.
$$
The reduced space
$(\C P^2\times \C P^2)_0$
has dimension 4, so we need to calculate the integral of $\kappa_0((\nu\otimes
1+ 1\otimes \nu)^2/2)$ over $(\C P^2\times \C P^2)_0$.
We instead find the pushforward of $[(\nu\otimes 1 + 1\otimes
\nu)^2/2][\eomega])$. Note that the term $\exp(\eomega)$ serves only as a
place-holder for each fixed point contribution. Thus we calculate the
pushforward by restricting to the fixed points. Note that
$$
a = \frac{(\nu\otimes 1+ 1\otimes
\nu)^2}{2}=\frac{1}{2}\left(\nu^2\otimes 1+2\nu\otimes \nu +
1\otimes \nu^2\right)
$$
 We begin with the
contribution at the fixed point $(S,S)$.

\begin{align*}
a|_{(S,S)} &=\frac{1}{2}\left((-u_1-u_2)^2\otimes
1+2(-u_1-u_2)\otimes (-u_1-u_2) + 1\otimes (-u_1-u_2)^2\right)\\
&= 2u_1^2 + 4u_1u_2+2u_2^2.
\end{align*}

We divide this restriction by the equivariant Euler class at $(S,S)$, which is described by the products of the
weights of the action. In this case, $e_{(S,S)}=u_1^2u_2^2$. Thus the contribution to the ABBV pushforward is
$$
\frac{(2u_1^2 + 4u_1u_2+2u_2^2)e^{-i(x_1+x_2)}}{u_1^2u_2^2} = 2\left(\frac{1}{u_2^2} +
2\frac{1}{u_1u_2} + \frac{1}{u_1^2}\right)e^{-i(x_1+x_2)}.
$$
We know by the linear story that we may Fourier transform this
rational function (assuming a polarization) and find the contribution
to the integral at 0 for this fixed point. We note that
$\frac{1}{u_2^2}e^{-i(x_1+x_2)}$ and $\frac{1}{u_1^2}e^{-i(x_1+x_2)}$ Fourier transform into
distributions supported on walls of the moment polytope, mainly on the
lines $u_1=-1$ and $u_2=-1$, respectively. Thus they do not contribute
to the integral. In contrast, the term $\frac{1}{u_1u_2}e^{-i(x_1+x_2)}$ contributes
``1'' on its support, which is the first quadrant out of the point
$(-1,-1)$ on $\t^*$. Thus the contribution of this class at $(S,S)$ to the integral
at 0 is 4.

Similarly we calculate the contributions of the other points. At both
$(S,E)$ and $(E,S)$ we obtain the restriction
\begin{align*}
a|_{(S,E)}&= \frac{1}{2}\left(1\otimes (2u_1-u_2)^2 + 2(-u_1-u_2)\otimes (2u_1-u_2)
+ (-u_1-u_2)^2\otimes 1\right)\\
&= \frac{1}{2}(u_1^2 - 4u_1u_2+4u_2^2).
\end{align*}
The equivariant Euler classes are:
$$
e_{(S,E)}= e_{(E,S)}= -u_1^2u_2(u_2-u_1)
$$
and so the contribution of each point to the integral of $\kappa_0(a)$ over $\C
P^2\times \C P^2)_0$ is
$$
\frac{ \frac{1}{2}(u_1^2 - 4u_1u_2+4u_2^2)}{-u_1^2u_2(u_2-u_1)}.
$$ We simplify and break up terms by partial fractions to obtain
\begin{align*}
&-e^{i(x_1-2x_2)}\frac{1}{2}\left(\frac{1}{u_2(u_2-u_1)} -
 \frac{4}{u_1(u_2-u_1)}+\frac{4u_2}{u_1^2(u_2-u_1)}\right)\\
&= -e^{i(x_1-2x_2)}\frac{1}{2}\left(\frac{1}{u_2(u_2-u_1)} -
 \frac{4}{u_1(u_2-u_1)}+ \frac{4}{u_1^2}+\frac{4}{u_1(u_2-u_1)}\right).
\end{align*}

The contribution at 0 is then $-\frac{1}{2}(1-4+4)= -\frac{1}{2}$ since $0$ is in the span of all the cones at
$(1,-2)$ indicated by the denominators except the term $\frac{4}{u^2}$. Thus the points $(S,E)$ and $(E,S)$
contribute a total of -1 to the symplectic volume. We now note that no other fixed points contribute to
integrals over $(\C P^2\times \C P^2)_0$ since the polarized cones out of other points do not contain 0. The
symplectic volume of a product of two copies of $\C P^2$ (with the given symplectic form) reduced by a diagonal
$T$ action is $4-1=3$.

\appendix

\section{The Equivariant K\"unneth Theorem}\label{se:Kunneth}

The following theorem is well-known to those who study equivariantly formal spaces, and hints toward other
proofs can be found, for example, in \cite{Hsiang}. The authors have been unable to find a satisfactory proof in
the literature.

\begin{theorem}[Equivariant K\"unneth Theorem]
Let $X$ and $Y$ be Hamiltonian $T$-spaces. Then in rational cohomology,
$$ H^*_T(\XxY) \iso H^*_T(X) \otimes_{H^*_T(pt)} H^*_T(Y) $$
as rings.
\end{theorem}

\begin{proof}
Let $M_T=M\times_TET$ be a bundle over $BT$ with fiber $M$, where $ET$ is a homotopically trivial space with a
free $T$ action and $BT$ is the quotient $ET/T$. This is sometimes called the {\em Borel space} associated to a
$T$-space $M$. The singular cohomology of $M_T$ is isomorphic to the equivariant Cartan cohomology $H_T^*(M)$
described in Section \ref{se:cartanmodel} \cite{AB:momentmap}. For a Hamiltonian $T$-space $M$, the Leray-Serre
spectral sequence associated to the bundle $M_T$ collapses at the $E_2$ term. Thus $H^*_T(M) \iso H^*(M) \otimes
H^*_T(pt) $ as modules over $H_T^*(pt)$ (and generally not as rings). In the case that $M=X\times Y$ is a
product of two Hamiltonian $T$-spaces considered under the diagonal $T$ action, this module isomorphism combined
with the ordinary K\"unneth theorem implies
$$
H_T^*(X\times Y)\iso H_T^*(X)\otimes_{H_T^*(pt)}H_T^*(Y)
$$ as modules.

To show the ring isomorphism, we find a map
$h: H^*_T(X) \otimes_{H^*_T} H^*_T(Y)\rightarrow H^*_T(\XxY)$.
Consider the commutative diagram
\begin{equation*}
\begin{CD}
(X\times Y)_T@>{\pi_1}>>X_T\\
@V{\pi_2}VV      @VV{\pi_X}V\\
Y_T@>{\pi_Y}>> BT
\end{CD}
\end{equation*}
We show that the surjective map $\pi^*_1\otimes\pi^*_2: H^*(X_T)\otimes
H^*(Y_T)\longrightarrow H^*((\XxY)_T)$ descends to a map
$$
h: H^*(X_T)\otimes_{H^*(BT)}
H^*(Y_T)\longrightarrow H^*((\XxY)_T).
$$
Let $u\in H^*(BT).$ Then $\pi_X^*(u)\in H^*(X_T)$ and
\begin{align*}
\pi^*_1\otimes\pi^*_2 (a\cdot \pi_X^*(u)\otimes b) &=
\pi_1^*(a)\pi_1^*\pi_X^*(u)\pi_2^*(b)\\
&= \pi_1^*(a)\pi_2^*\pi_Y^*(u)\pi_2^*(b)\\
&=\pi^*_1\otimes\pi^*_2 (a \otimes \pi_Y^*(u)\cdot b),
\end{align*}
which proves that $h$ exists and is an isomorphism of rings.
\end{proof}

\end{document}